\documentclass{amsproc}
\usepackage{graphicx}
\usepackage{amssymb}
\usepackage{epstopdf}
\usepackage{amsmath,amsthm,amscd,amssymb}
\usepackage{latexsym}
\usepackage[colorlinks,citecolor=red,pagebackref,hypertexnames=false]{hyperref}
\usepackage{geometry}                
\numberwithin{equation}{section}

\theoremstyle{plain}
\newtheorem{theorem}{Theorem}[section]

\theoremstyle{definition}

\newtheorem{case[theorem]}{Case}

\theoremstyle{remark}

\def\od{\mathbb O(d)}

\def\R{\mathbb R}

\numberwithin{equation}{section}

\begin{document}

\title{An improved dimensional threshold for the angle problem}

\author{Alex Iosevich and Eyvindur A. Palsson}

\date{today}

\email{iosevich@math.rochester.edu}
\email{palsson@vt.edu}

\address{Department of Mathematics, University of Rochester, Rochester, NY 14627}
\address{Department of Mathematics, Virginia Tech, Blacksburg, VA 24061}

\thanks{The work of the first listed author was partially supported by the NSA Grant H98230-15-1-0319. The work of the second listed author was supported in part by Simons Foundation Grant \#360560.}

\begin{abstract} 
The Falconer distinct distance problem asks for a compact set $E\subset\mathbb{R}^d$ how large its Hausdorff dimension needs to be to ensure that the Lebesgue measure of its distance set is positive. In this paper we consider the analogous question for the set of angles. We show that if the Hausdorff dimension of $E$ is strictly bigger than $\frac{d}{2}$ then the Lebesgue measure of the angles set is positive. In the plane this result was previously established by Harangi et al \cite{HKKMMMS13}. In higher dimensions, our exponent improves the $\frac{d+1}{2}$ threshold previously obtain by the authors of this paper and Mihalis Mourgoglou \cite{IMP16}. We do not know what the right dimensional threshold should be in higher dimensions. 
\end{abstract} 

\maketitle


\section{Introduction}

\vskip.125in 

One of the most important and far reaching problems in modern geometric measure theory is the Falconer distance problem, which asks: How large does the Hausdorff dimension $s$ of a compact set $E \subset {\mathbb R}^d$, $d \ge 2$, need to be to ensure that the \emph{distance set} of $E$, $\Delta(E):=\{|x-y|: x,y \in E \}\subset \R$, has positive Lebesgue measure? Falconer proved that $s>\frac{d}{2}$ is necessary, up to the endpoint, and conjectured that it is also sufficient \cite{Falc86}. In his original paper Falconer obtained the threshold $\frac{d}{2}+\frac{1}{2}$ through an incidence theorem which heuristically says that no distance appears particularly often. This can be interpreted as an $L^{\infty}$ style argument. Wolff improved this result in the plane to $\frac{4}{3}$ by showing that the probability that any particular distance would appear had an $L^2$ density with respect to the Lebesgue measure \cite{W99}. His result was then extended to higher dimension by Erdo\~{g}an, who obtained the threshold $\frac{d}{2}+\frac{1}{3}$ \cite{Erd05}. This was improved in dimensions $3$ and higher by 
Du, Guth, Ou, Wang, Wilson and Zhang \cite{DGOWWZ}, who obtained the threshold $1.8$ in dimension $d=3$ and $\frac{d}{2}+\frac{1}{4}+\frac{d+1}{4(2d+1)(d-1)}$ in dimensions $d\geq 4$, and then further improved in dimensions $4$ and higher by Du and Zhang \cite{DZ} who obtained $\frac{d}{2}+\frac{1}{4}+\frac{1}{8d-4}$ for $d\geq 2$ which additionally matches the best known results in dimensions $2$ and $3$. These improvements were obtained through improved Fourier restriction estimates, while still using the same setup as Wolff and Erdo\u{g}an. Most recently Guth, Iosevich, Ou and Wang \cite{GIOW} obtained the threshold $\frac{5}{4}$ in dimension $d=2$ where understanding the sharpness examples for the Fourier restriction estimates played a key role. For related results on the dimension of the distance set see Orponen \cite{O17}, Shmerkin \cite{Shm17} and Keleti-Shmerkin \cite{KS18}. 

Similar questions, as Falconer proposed for distance, can be asked for more general point configurations. Greenleaf and the first named author initiated the study of more general configurations when they considered triangles in the plane \cite{GI12}. This has been expanded to simplices \cite{GILP} and even more general configurations \cite{GGIP12}. A particularly relevant paper is that of the two authors with Greenleaf and Liu \cite{GILP} where they took a novel group-theoretic approach to the Wolff and Erdo\u{g}an approach to simplices, which shed some new light even on the case of distances. This motivated the result and approach in this paper.

A particular point configuration of interest is that of angles. Let $E\subseteq\mathbb{R}^d$ be compact and define the angles set
$$\mathcal{A}(E)=\left\lbrace \theta(x,y,z): x,y,z \in E \text{ distinct} \right\rbrace$$
where $\theta(x,y,z)$ denotes the angle formed by the points $x,y,z$ centered at $y$ sitting in $[0,2\pi]$. A Falconer type question is how large does $dim_{{\mathcal H}}(E)$ need to be to ensure that the Lebesgue measure of the angles set is positive, $\mathcal{L}(\mathcal{A}(E))>0$. This configuration was studied by Mourgoglou and the two named authors in \cite{IMP16} where they obtained the threshold $\frac{d+1}{2}$ through an incidence theorem. Extending a construction by Apfelbaum and Sharir in the discrete setting to higher dimensions and the continuous setting it was also shown that the incidence theorem was sharp for the angle $\frac{\pi}{2}$. Harangi et al \cite{HKKMMMS13} and M\'{a}th\'{e} \cite{Mathe17} studied how large a dimension guarantees a given angle. A corollary of their results is the threshold $d-1$, which in particular yields the threshold $1$ in the plane which is sharp, since if $E$ is a line in the plane then $\mathcal{A}(E)$ only contains two angles, $0$ and $\pi$, so $\mathcal{L}(\mathcal{A}(E))=0$.

The first result in this paper takes the approach of Wolff and contains the essential features of our method, where we have adapted the group theoretic approach of \cite{GILP} to angles. Challenges involve more complicated geometry and a distinct role of the mid-point of the angle, which is in contrast to what happened in \cite{GILP} where all points had a similar role.

\begin{theorem} \label{method}
Let $E$ be a compact set in $ {\mathbb R}^d, \, d \ge 2$, and $\mu$ a finite, nonnegative Frostman measure supported on $E$. For $a\in\R^+$ and $g\in\od $, the orthogonal group on $\R^d$, define a measure  $\nu_{a,g}$, supported  on $E-agE$, by the relation 
\begin{equation} \label{rotatedmeasure} \int_{\R^d} f(z)\, d\nu_{a,g}(z):=\int_E \int_E f(u-agv)\, d\mu(u)\, d\mu(v),\, f\in C_0(\R^d). \end{equation} 
Define  also a measure $\nu$ on $\mathcal{A}(E)\subset\mathbb{R}$ by 
\begin{equation} \label{anglemeasure} \int f(t)\, d\nu(t)=\iiint f\left( \theta(x,y,z) \right) d\mu(x) d\mu(y) d\mu(z), \end{equation} where $t$ is which is the push-forward of $\mu\times\mu\times\mu$ under the map $(x,y,z) \mapsto \theta(x,y,z)$.

\vskip.125in 

Then, if there exists a compact interval $I\subset  {\mathbb R}^{+}$ such that $\nu_{a,g}$ is absolutely continuous for $\hbox{a.e.} \, (a,g)\in I\times\od$, with density also denoted $\nu_{a,g}$, and
\begin{equation} \label{knorm} \int_I \int_{\od} \int_{\R^d} \nu_{a,g}^{2}(x)\, dx\,dg\, \frac{da}{a}<\infty, \end{equation}
where $dg$ is Haar measure on $\od$, 
then the measure $\nu$ in (\ref{anglemeasure}) has an $L^2$ density and    ${\mathcal L}(\mathcal{A}(E))>0$. 

\end{theorem}

Here we remind the reader that $\mu$ is a Frostman measure supported on $E$ (see e.g. \cite{W03}, Chapter 8) if for any $\epsilon>0$ there exists $C_{\epsilon}>0$ such that if $B_{\delta}$ is a ball of radius $\delta$, then $\mu(B_{\delta}) \leq C_{\epsilon}\delta^{s-\epsilon}$ where $s=dim_{{\mathcal H}}(E)$.

As an application of Theorem \ref{method} we match the sharp threshold $1$ in the plane and improve the dimensional threshold $\frac{d+1}{2}$ of \cite{IMP16} to $\frac{d}{2}$ in higher dimensions.

\begin{theorem} \label{mainangle} Let $E \subset \mathbb{R}^d$ be compact, $d \geq 2$. Suppose that $dim_{{\mathcal H}}(E)> \frac{d}{2}$. Then $\mathcal{L}(\mathcal{A}(E))>0$.
\end{theorem}

\vskip.25in

\section{Proof of Theorem \ref{method}}

We adapt the techniques of \cite{GILP} to the setting of angles. Define a measure $d\nu$ on $\R$, with support in $\mathcal{A}(E)$, as in (\ref{anglemeasure}) above. We will show that  to prove Theorem \ref{method} it  suffices to obtain an upper bound on the $L^2$ norm of the density, i.e., the Radon-Nikodym derivative of  $d\nu$, which we denote by $\nu(t)$. We start by showing that
\begin{equation}\label{thickened}
\int \nu^2(t)\, dt\le c_{d}\cdot \liminf_{\epsilon \to 0} \epsilon^{-1} \mu^{6} \big\{ x_1,y_1,z_1,x_2,y_2,z_2\in\mathbb{R}^d: \left|\theta(x_1,y_1,z_1) - \theta(x_2,y_2,z_2)\right|\lesssim \epsilon \big\},
\end{equation}
where $\mu^{6}$ denotes $\mu \times \mu \times \mu \times \mu \times \mu \times \mu$, with the proof showing that if the RHS of (\ref{thickened}) is finite,
then in fact $d\nu$ is absolutely continuous with respect to Lebesgue measure $dt$, with density $\nu(t)\in L^2$.

Let $\phi\in C_0^\infty(\R),\, \phi\ge 0,\, supp(\phi)\subset \{|t|\le 1\},\, \int \phi \, dt=1$, and $\phi_\epsilon(\cdot)=\epsilon^{-1}\phi(\epsilon^{-1}\cdot),\, 0<\epsilon<\infty$, the resulting approximate identity. 
Setting $\nu_\epsilon=\phi_\epsilon * d\nu\in C_0^\infty$, one has $d\nu=wk^*\!-\!\lim_{\epsilon\to 0} \nu_\epsilon$, and
(\ref{thickened}) will follow if one shows that $\liminf_{\epsilon\to 0} ||\nu_\epsilon||_{L^2}^2=C<\infty$.

\vskip.125in 

Now, 
$$\nu_\epsilon(t) = \iiint \phi_\epsilon\left( \theta(x,y,z) - t \right)\, d\mu(x) d\mu(y) d\mu(z).$$
Due to the nonnegativity of $\phi_\epsilon$ and $d\mu$, this is dominated by
$$\iiint\, \epsilon^{-1}\chi\left\{\left|\theta(x,y,z) - t\right|<\epsilon\right\}\, d\mu(x) d\mu(y) d\mu(z),$$
where $\chi(A)$ denotes the characteristic function of a set $A$, and thus

\begin{eqnarray}\label{squared}
||\nu_\epsilon||_{L^2}^2 & \lesssim & \idotsint \epsilon^{-1}\chi\left\{\left|\theta(x_1,y_1,z_1) - t\right|<\epsilon\right\}\, \epsilon^{-1}\chi\left\{\left|\theta(x_2,y_2,z_2) - t\right|<\epsilon\right\} \\
& & \nonumber \qquad \qquad \qquad \qquad \qquad \qquad d\mu(x_1)\, d\mu(y_1)\, d\mu(z_1)\, d\mu(x_2)\, d\mu(y_2)\, d\mu(z_2) \,dt.
\end{eqnarray}
Now, by the triangle inequality, one has
\begin{equation*}
\chi\left\{\left|\theta(x_1,y_1,z_1) - t\right|<\epsilon\right\}\cdot\chi\left\{\left|\theta(x_2,y_2,z_2) - t\right|<\epsilon\right\} \leq
\chi\left\{\left|\theta(x_1,y_1,z_1)-\theta(x_2,y_2,z_2)\right|<2\epsilon\right\},
\end{equation*}
and thus, integrating out $dt$, the RHS of (\ref{squared}) is
$$\lesssim \epsilon^{-1}\idotsint \chi\left\{\left|\theta(x_1,y_1,z_1)-\theta(x_2,y_2,z_2)\right|<2\epsilon\right\} d\mu(x_1)\, d\mu(y_1)\, d\mu(z_1)\, d\mu(x_2)\, d\mu(y_2)\, d\mu(z_2).$$
Taking the $\liminf$ as $\epsilon\to 0$ yields the RHS of (\ref{thickened}).

\vskip.125in

Fix $x_1,y_1,z_1,x_2,y_2,z_2$ as in (\ref{squared}). There exists a rotation $g\in\od$ such that
$$ \frac{x_1-y_1}{|x_1-y_1|} = g \frac{x_2-y_2}{|x_2-y_2|} .$$
The vectors $\frac{x_1-y_1}{|x_1-y_1|} = g \frac{x_2-y_2}{|x_2-y_2|}$ and $\frac{z_1-y_1}{|z_1-y_1|}$ span a plane. There exists a rotation $g^{*}\in\text{Stab}\left( \frac{x_1-y_1}{|x_1-y_1|} \right)$ such that $g^{*}g\frac{z_2-y_2}{|z_2-y_2|}$ sits in that plane, note this step is unnecessary if in the case $d=2$ and is only relevant if $d\geq 3$. Here $\text{Stab}(u)$ is the stabilizer of the vector $u$, that is, the set of rotations from $\od$ that leave the vector $u$ unchanged. Note that we can identify the subgroup $\text{Stab}\left( \frac{x_1-y_1}{|x_1-y_1|}\right)$ with $\mathbb{O}(d-1)$. Note that the angle between the vectors $\frac{x_1-y_1}{|x_1-y_1|}$ and $\frac{z_1-y_1}{|z_1-y_1|}$ is $\theta(x_1,y_1,z_1)$ and the angle between $g \frac{x_2-y_2}{|x_2-y_2|} = \frac{x_1-y_1}{|x_1-y_1|}$ and $g^{*}g\frac{z_2-y_2}{|z_2-y_2|}$ is $\theta(x_2,y_2,z_2)$ because rotations preserve angles and $g^{*}\in\text{Stab}\left( \frac{x_1-y_1}{|x_1-y_1|} \right)$. Thus we see that the angle between $\frac{z_1-y_1}{|z_1-y_1|}$ and $g^{*}g\frac{z_2-y_2}{|z_2-y_2|}$ is $\omega := \theta(x_1,y_1,z_1) - \theta(x_2,y_2,z_2)$ and we know $|\omega| \lesssim \epsilon$. Note $\frac{z_1-y_1}{|z_1-y_1|}$ and $g^{*}g\frac{z_2-y_2}{|z_2-y_2|}$ form an isosceles triangle with the two equal sides having length $1$ and the angle between them being $\omega$. By simple trigonometry we see that the third side, which corresponds to the distance between the two vectors, has length $|2\sin(\omega/2)|$ which then immediately yields
$$ \left| \frac{z_1-y_1}{|z_1-y_1|} - g^{*}g\frac{z_2-y_2}{|z_2-y_2|} \right| \lesssim \epsilon .$$
Finally we can write
$$ \frac{x_1-y_1}{|x_1-y_1|} = g \frac{x_2-y_2}{|x_2-y_2|} \, \text{ as } \, x_1-y_1 = ag(x_2 - y_2) $$
for a scalar $a\in {\mathbb R}^{+}$ and similarly we can write
$$ \left| \frac{z_1-y_1}{|z_1-y_1|} - g^{*}g\frac{z_2-y_2}{|z_2-y_2|} \right| \lesssim \epsilon \, \text{ as } \, |(z_1-y_1) - bg^{*}g(z_2-y_2) |\lesssim \epsilon $$
for a scalar $b\in {\mathbb R}^{+}$. Observe that due to pigeon holing, where the set $E$ is broken up into three separated sets with positive $\mu$ measure, the different $a$'s and $b$'s can be bounded below with an absolute constant that is strictly positive and likewise they can be bounded above with the reciprocal of the lower bound times twice the diameter of the set $E$. Thus the different $a$'s and $b$'s can be taken from a compact interval $I$ that is away from $0$.

Now for each $z_2 - y_2$ take a cover of $\od/\text{Stab}(z_2 - y_2)$ by balls of radius $\epsilon$ with respect to the induced Riemannian metric with finite overlap. Since the dimension of $\od/\text{Stab}(z_2 - y_2)$ is that of $\od / \mathbb{O}(d-1)$, namely $d-1$, we need $N(\epsilon) \sim C \epsilon^{-(d-1)}$ balls to cover it. Choose sample points from $\od/\text{Stab}(z_2 - y_2) \times I \times I$, $(\tilde{g}_m, a_m, b_m)$, $1\leq m \leq N(\epsilon)$, where $\tilde{g}_m = \tilde{g}_m(z_2 - y_2)$ is taken one in each of the balls. We thus see that the set
$$ \lbrace x_1,y_1,z_1,x_2,y_2,z_2 : \left|\theta(x_1,y_1,z_1)-\theta(x_2,y_2,z_2)\right| \lesssim \epsilon \rbrace $$
is contained in
\begin{multline*}
\bigcup\limits_{m=1}^{N(\epsilon)} \biggl\lbrace x_1,y_1,z_1,x_2,y_2,z_2 : |(x_1 - y_1) - a_m\tilde{g}_mh(x_2-y_2)| \lesssim \epsilon \\
\text{ and } |(z_1-y_1) - b_m g^{*}\tilde{g}_m h(z_2-y_2) |\lesssim \epsilon, h\in\text{Stab}(z_2 - y_2) \biggr\rbrace .
\end{multline*}

Since this holds for any choice of sample points $(\tilde{g}_m, a_m, b_m)$, we can pick these points such that they minimize (up to a factor of $1/2$, say) the quantity
\begin{multline*}
\mu^{6} \biggl\lbrace x_1,y_1,z_1,x_2,y_2,z_2 : |(x_1 - y_1) - a_m\tilde{g}_mh(x_2-y_2)| \lesssim \epsilon \\
\text{ and } |(z_1-y_1) - b_m g^{*}\tilde{g}_m h(z_2-y_2) |\lesssim \epsilon, h\in\text{Stab}(z_2 - y_2) \biggr\rbrace .
\end{multline*}

Now consider the $N(\epsilon)$ preimages, under the natural projection from $\od$, of the balls used to cover $\od / \text{Stab}(z_2 - y_2)$; we can label  these $\epsilon$-tubular neighborhoods of the preimages of the sample points $\tilde{g}_m$ as $T_1^{\epsilon},\ldots,T_{N(\epsilon)}^{\epsilon}$.  
Since $dim(\od / \text{Stab}(z_2 - y_2))=d-1$, each  $T^\epsilon_m$ has volume  \linebreak$\sim \epsilon^{d-1}$. The inf over a set is less than or equal to the average over the set, so we obtain that
\begin{multline*}
\mu^{6} \biggl\{x_1,y_1,z_1,x_2,y_2,z_2 :  |(x_1 - y_1) - a\tilde{g}_mh(x_2-y_2)| \lesssim \epsilon \\
\text{ and } |(z_1-y_1) - b g^{*}\tilde{g}_m h(z_2-y_2) |\lesssim \epsilon, h\in\text{Stab}(z_2 - y_2) \biggr\}
\end{multline*}
is bounded above, up to constants that depend on the length of the interval $I$, by
\begin{multline*}
\int\limits_{I} \int\limits_{I}\frac{1}{\epsilon^{d-1}}\int\limits_{T_m^{\epsilon}} \mu^{6} \biggl\{x_1,y_1,z_1,x_2,y_2,z_2 :  |(x_1 - y_1) - a g(x_2-y_2)| \lesssim \epsilon  \text{ and } |(z_1-y_1) - b g^{*}g(z_2-y_2) |\lesssim \epsilon \biggr\} \, dg \, \frac{da}{a} \frac{db}{b}
\end{multline*}
which allows us to bound the expression within the $\liminf$ on the RHS of (\ref{thickened}) above by
\begin{multline*}
\int\limits_{I} \int\limits_{I} \, \epsilon^{-d} \sum\limits_{m=1}^{N(\epsilon)} \, \int\limits_{T_m^{\epsilon}} \mu^{6} \biggl\{x_1,y_1,z_1,x_2,y_2,z_2 :  |(x_1 - y_1) - a g(x_2-y_2)| \lesssim \epsilon \\
\text{ and } |(z_1-y_1) - b g^{*}g(z_2-y_2) |\lesssim \epsilon \biggr\} \, dg \, \frac{da}{a} \frac{db}{b}
\end{multline*}
We note that the usage of the Haar measure is done for convenience for later calculations and is up to constants equivalent to using the Lebesgue measure as our compact interval $I$ is bounded away from $0$.
Since the cover has finite overlap, this in turn can be bounded above, up to a  constant $c_{d}$,  by
$$ \epsilon^{-d} \int\limits_{I} \int\limits_{I} \int\limits_{\od} \mu^{6} \biggl\{x_1,y_1,z_1,x_2,y_2,z_2 :  |(x_1 - y_1) - a g(x_2-y_2)| \lesssim \epsilon
\text{ and } |(z_1-y_1) - b g^{*}g(z_2-y_2) |\lesssim \epsilon \biggr\} \, dg \, \frac{da}{a} \frac{db}{b}.$$
Taking the limit as $\epsilon\rightarrow 0^{+}$, we obtain a constant multiple of the expression
$$ \int\limits_{I}\int\limits_{I}\int\limits_{\od} \mu^{6} \biggl\{x_1,y_1,z_1,x_2,y_2,z_2 :  (x_1 - y_1) = a g(x_2-y_2) \text{ and } (z_1-y_1) = b g^{*}g(z_2-y_2) \biggr\} \, dg \, \frac{da}{a} \, \frac{db}{b}. $$
Now write this as
\begin{multline*}
\biggl( \biggl(\int\limits_{I}\int\limits_{I}\int\limits_{\od}\int\limits_{\mathbb{R}^6} \chi\left\lbrace (x_1 - y_1) = a g(x_2-y_2) \text{ and } (z_1-y_1) = b g^{*}g(z_2-y_2) \right\rbrace \\
d\mu(x_1)\, d\mu(y_1)\, d\mu(z_1)\, d\mu(x_2)\, d\mu(y_2)\, d\mu(z_2) \, dg \, \frac{da}{a} \, \frac{db}{b} \biggr)^{1/2} \biggr)^{2} \\
\end{multline*}
and then enlarge the area of integration on the one hand by dropping $(z_1-y_1) = b g^{*}g(z_2-y_2)$ and on the other by dropping $(x_1 - y_1) = a g(x_2-y_2)$ and thus obtain the upper bound
\begin{multline*}
\biggl(\int\limits_{I}\int\limits_{I}\int\limits_{\od}\int\limits_{\mathbb{R}^{6d}} \chi\left\lbrace (x_1 - y_1) = a g(x_2-y_2) \right\rbrace
d\mu(x_1)\, d\mu(y_1)\, d\mu(z_1)\, d\mu(x_2)\, d\mu(y_2)\, d\mu(z_2) \, dg \, \frac{da}{a} \, \frac{db}{b} \biggr)^{1/2} \\
\cdot \biggl(\int\limits_{I}\int\limits_{I}\int\limits_{\od}\int\limits_{\mathbb{R}^{6d}} \chi\left\lbrace (z_1-y_1) = b g^{*}g(z_2-y_2) \right\rbrace
d\mu(x_1)\, d\mu(y_1)\, d\mu(z_1)\, d\mu(x_2)\, d\mu(y_2)\, d\mu(z_2) \, dg \, \frac{da}{a} \, \frac{db}{b} \biggr)^{1/2}
\end{multline*}
and finally integrate through with the free variables, which is possible because $I$ is finite and away from $0$ and $\mu$ is a probability measure, and obtain an upper bound, up to constants, of
\begin{multline*}
\biggl(\int\limits_{I}\int\limits_{\od} \int\limits_{\mathbb{R}^{4d}} \chi\left\lbrace x_1 - a g x_2 = y_1 - a g y_2 \right\rbrace
d\mu(x_1)\, d\mu(y_1)\, d\mu(x_2)\, d\mu(y_2) \, dg \, \frac{da}{a} \biggr)^{1/2} \\
\cdot \biggl(\int\limits_{I}\int\limits_{\od}  \int\limits_{\mathbb{R}^{4d}} \chi\left\lbrace z_1 -  b g z_2 = y_1 - b g y_2 \right\rbrace
d\mu(y_1)\, d\mu(z_1)\, d\mu(y_2)\, d\mu(z_2) \, dg \, \frac{db}{b} \biggr)^{1/2}
\end{multline*}
where we also used that $g^{*}g$ ranges through all of $\od$ as $g$ ranges through all of $\od$. Finally using the definition of $\nu_{a,g}$ from (\ref{anglemeasure}) we obtain the following bound
\begin{equation}
\int \nu^2(t)\, dt \lesssim \int\limits_{I}\int\limits_{\od} \int\limits_{\mathbb{R}^{d}} \nu_{a,g}^{2}(x) \, dx \, dg \, \frac{da}{a} ,
\end{equation}
which is (\ref{knorm}) as we wanted to prove. This completes the proof of Theorem \ref{method}.

\vskip.25in

\section{Proof of Theorem \ref{mainangle}} \label{sec mainproof}

\vskip.125in 

The matters have been reduced in the introduction to the estimation of (\ref{knorm}). Let $\psi$ be a smooth cutoff function supported in $\left\{\xi \in {\mathbb R}^d: \frac{1}{2} \leq |\xi| \leq 4 \right\}$ and identically equal to $1$ in $\left\{\xi \in {\mathbb R}^d: 1 \leq |\xi| \leq 2  \right\}$. Let $\nu_{a,g,j}$ denote the $j$th Littlewood-Paley piece of $\nu_{a,g}$, defined by the relation $\widehat{\nu}_{a,g,j}(\xi)=\widehat{\nu}_{a,g}(\xi) \psi(2^{-j}\xi)$. Since $\nu_{a,g}$ is compactly supported, we may assume that $j\ge 0$. Using the Littlewood-Paley decomposition of $\nu_{a,g}$, the inner part of the integral in (\ref{knorm}) equals
$$ \int \sum_{j_1, j_2} \nu_{a,g,j_1}(x) \nu_{a,g,j_2}(x) \,  dx ,$$
which by Plancherel equates to
$$  \sum_{j_1, j_2} \int \widehat{\nu}_{a,g,j_1}(\xi) \widehat{\nu}_{a,g,j_2}(\xi) \, d\xi $$
so the sum vanishes if $|j_1 - j_2| > 2$. Thus it suffices to consider the case $j_1 = j_2$ and study
$$ \sum_{j} \int \nu_{a,g,j}^2(x) \, dx .$$
Using the definition of $\nu_{t,g}$ with $f(z)=e^{-2 \pi i z \cdot \xi}$, we obtain
$$ \widehat{\nu}_{t,g}(\xi)=\widehat{\mu}(\xi) \overline{\widehat{\mu}(t g \xi)},$$
which means that, via Plancherel,
$$ \int\limits_{I}\int\limits_{\od} \int\limits_{\mathbb{R}^{d}} \nu_{a,g,j}^{2}(x) \, dx \, dg \, \frac{da}{a} = \int_{2^{j}\leq |\xi| \leq 2^{j+1}}{|\widehat{\mu}(\xi)|}^2 \left\{ \int\limits_{I}\int\limits_{\od} {|\widehat{\mu}(a g \xi)|}^2 dg \, \frac{da}{a} \right\} d\xi. $$
Switching to polar coordinates and using the action of the orthogonal group on the sphere, one quickly sees that this quantity equals a constant multiple of 
\begin{equation*}
\int_{2^j}^{2^{j+1}} {\left( \int\limits_{I}\int\limits_{S^{d-1}} {|\widehat{\mu}(a t \omega)|}^2 d\omega \, \frac{da}{a} \right)}^2 t^{d-1} dt.
\end{equation*}
Now observe that
\begin{eqnarray*}
 \int_{I} \int_{S^{d-1}} {|\widehat{\mu}(a t \omega)|}^2\, d\omega \, \frac{da}{a} &=&\int_{tI} \int_{S^{d-1}} {|\widehat{\mu}(a \omega)|}^2\, d\omega\,  \frac{da}{a} \\
& \lesssim &  {t}^{-d}  \int_{tI} \int_{S^{d-1}} {|\widehat{\mu}(a \omega)|}^2 \, d\omega\, a^{d-1}\,  da\\
&\lesssim &{t}^{-d} \int_{|x|\in tI} {|\widehat{\mu}(x)|}^2\,  dx \lesssim {t}^{-s}
\end{eqnarray*}
with implicit constants depending on $I$, where we use that the last expression is an energy integral and $\mu$ is a Frostman measure supported on the set $E$ with $dim_{{\mathcal H}}(E)>s$. Plugging this into the integral above yields
$$ \int_{2^j}^{2^{j+1}} {\left( \int\limits_{I}\int\limits_{S^{d-1}} {|\widehat{\mu}(a t \omega)|}^2 d\omega \, \frac{da}{a} \right)}^2 t^{d-1} dt \lesssim (2^{j})^{d-2s} $$
which shows that the geometric series converges if $s>\frac{d}{2}$, which implies $dim_{{\mathcal H}}(E)>\frac{d}{2}$. For related calculations see \cite{Mat87}.

\bigskip


\vskip.25in 


 

\begin{thebibliography}{8}

\bibitem{DGOWWZ}
X. Du, L. Guth, Y. Ou, H. Wang, B. Wilson, R. Zhang,
\emph{Weighted restriction estimates and application to Falconer distance set problem}, (2018), arXiv:1802.10186.

\bibitem{DZ}
X. Du, R. Zhang,
\emph{Sharp $L^2$ estimate of Schr\"{o}dinger maximal function in higher dimensions}, (2018), arXiv:1805.02775.

\bibitem{Erd05} B. Erdo\~{g}an {\it A bilinear Fourier extension theorem and applications to the distance set problem}, Int. Math. Res. Not.  (2006).

\bibitem{Falc86} K. J. Falconer, {\it On the Hausdorff dimensions of distance sets}, Mathematika \textbf{32} (1986), 206-212.

\bibitem{GIOW}
L. Guth, A. Iosevich, Y. Ou, H. Wang,
\emph{On Falconer's distance set problem in the plane}, (2018), arXiv:1808.09346.

\bibitem{GGIP12} L. Grafakos, A. Greenleaf, A. Iosevich and E. Palsson, {\it Multilinear generalized Radon transforms and point configurations}, Forum Mathematicum, \textbf{27} (2015), 2323-2360. 

\bibitem{GI12} A. Greenleaf, A. Iosevich, {\it On three point configurations determined by subsets of the Euclidean plane, the associated bilinear operator and applications to discrete geometry}, Analysis and PDE, \!\textbf{5}-2\! (2012), 397-409. 

\bibitem{GILP} A. Greenleaf, A. Iosevich, B. Liu, E. A. Palsson, \emph{A group-theoretic viewpoint on Erd\H{o}s-Falconer problems and the Mattila integral}, Revista Matematica Iberoamericana, \textbf{31} (2015), no. 3, 799-810.

\bibitem{HKKMMMS13} V. Harangi, T. Keleti, G. Kiss, P. Maga, A. Mathe, P. Mattila, B. Strenner, {\it How large dimension guarantees a given angle?} (English summary) Monatsh. Math. \textbf{171} (2013), no. 2, 169-187. 

\bibitem{IMP16} A. Iosevich, M. Mourgoglou, E. A. Palsson, \emph{On angles determined by fractal subsets of the Euclidean space via Sobolev bounds for bi-linear operators}, Mathematical Research Letters, \textbf{23} (2016), 1737-1759.

\bibitem{KS18} T. Keleti and P. Shmerkin, {\it New bounds on the dimensions of planar distances sets}, (2018) arXiv:1801.08745. 

\bibitem{Mathe17} A. Mathe, {\it Sets of large dimension not containing polynomial configurations}, Adv. Math. 316 (2017), 691-709.

\bibitem{Mat87} P. Mattila {\it Spherical averages of Fourier transforms of measures with finite energy: dimensions of intersections and distance sets}, Mathematika, \textbf{34} (1987),  207-228.

\bibitem{O17} T. Orponen, {\it On the distance sets of Ahlfors-David regular sets}, Advances in Mathematics 307 (2017), 1029-1045.

\bibitem{Shm17} P. Shmerkin, {\it On distance sets, box-counting and Ahlfors regular sets}, Discrete Anal. (2017), No. 9, 22 pp.

\bibitem{W99} T. Wolff, {\it Decay of circular means of Fourier transforms of measures}, Int. Math. Res. Not.  \textbf{10} (1999) 547--567.

\bibitem{W03} T. Wolff, \emph{Lectures on harmonic analysis}, I. Laba and C. Shubin, eds. University Lecture Series, \textbf{29}. Amer. Math. Soc., Providence, RI, 2003.



\end{thebibliography}
\end{document}